\documentclass[12pt]{article}
\usepackage{mathptmx}
\usepackage{amsfonts}
\newtheorem{thm}{Theorem}[section]
\newtheorem{lem}{Lemma}[section]
\newtheorem{rmk}{Remark}[section]
\newtheorem{prop}{Proposition}[section]

\newtheorem{pf}{Proof}

\title{Approximation by generalized bivariate Kantorovich sampling type series}
\author{A. Sathish Kumar\thanks{Department of Mathematics, Visvesvaraya National Institute of Technology, Nagpur, Nagpur-440010, India. \newline E-mail: mathsatish9@gmail.com}
 \and
 P. Devaraj\thanks{Department of Mathematics, Indian Institute of Science Education and Research, Thiruvananthapuram \newline E-mail: devarajp@iisertvm.ac.in}
  }
\date{}

\begin{document}
\maketitle
\bibliographystyle{plain}
\abstract{The purpose of this paper is to construct a bivariate generalization of new family of Kantorovich type sampling operators $(K_w^{\varphi}f)_{w>0}.$ First, we give the pointwise convergence theorem and a Voronovskaja type theorem for these Kantorovich generalized sampling series. Further, we obtain the degree of approximation by means of modulus of continuity and quantitative version of Voronovskaja type theorem for the family $(K_w^{\varphi}f)_{w>0}.$ Finally, we give some examples of kernels such as box spline kernels and Bochner-Riesz kernel to which the theory can be applied.

\endabstract



\section{Introduction}
In recent years, theory of generalized sampling series is an attractive topic in approximation theory due to its wide range of applications, especially in signal and image processing. The theory of generalized sampling series was first initiated by P. L. Butzer and his school \cite{Stens1} and \cite{Stens2}. After this many researchers have studied in this direction and obtained the convergence properties of these operators (cf. \cite{a7,a5,a3,a1} etc.).

Also, the Kantorovich type generalizations of the generalized sampling operators were introduced by P. L. Butzer and his school (see e.g. \cite{Stens1,PBL3,Don,GV1} etc). In \cite{BVBS}, the authors introduced the sampling Kantorovich operators and studied their rate of convergence in the general settings of Orlicz spaces. Danilo and Vinti \cite{Dani} obtained the rate of approximation for the family of sampling Kantorovich operators in the uniform norm for uniformly continuous bounded functions belonging to Lipschitz classes and for functions in Orlicz spaces. Also, the nonlinear version of sampling Kantorovich operators have been studied in \cite{Dan} and
\cite{Zam}.

Let $\varphi:{\mathbb{R}}^2\rightarrow {\mathbb{R}}$ be a suitable kernel function.  Then, the two dimensional generalized sampling series of a function  $f:{\mathbb{R}}^2\rightarrow {\mathbb{R}}$ is defined by
\begin{eqnarray*}\label{e1}
(T_w^{\varphi}f)(x,y)&=&\sum_{k=-\infty}^{\infty}\sum_{j=-\infty}^{\infty}\varphi\bigg(w\bigg(x-\frac{k}{w}\bigg), w\bigg(y-\frac{j}{w}\bigg)\bigg)f\bigg(\frac{k}{w},\frac{j}{w}\bigg),
\end{eqnarray*}
where $w\in \mathbb{N}$ and $(x,y)\in {\mathbb{R}}^2.$ These operators have great importance in the development of models for signal recovery. These type of operators have been studied by many authors (cf. \cite{bv1,PLB2,PLB1,Fis} etc.).

Altomare and Leonessa \cite{Alt} considered a new sequence of positive linear operators acting on the space of Lebesgue-integrable functions on the unit interval. Such operators include the Kantorovich operators as a particular case. Later, in order to obtain an approximation process for spaces of locally integrable functions on unbounded intervals, Altomare et. al. introduced and studied the generalized  Sz\'{a}sz-Mirakjan-Kantorovich operators in \cite{MMC}. Also in \cite{Cap}, the authors obtained some qualitative properties and an asymptotic formula for such a sequence of operators.

Motivated by the above works, we consider the bivariate generalized Kantorovich sampling series. Let $\{a_k\}_{k\in {\mathbb{Z}}}, \{b_k\}_{k\in {\mathbb{Z}}},\{c_j\}_{j\in {\mathbb{Z}}}$ and $\{d_j\}_{j\in {\mathbb{Z}}}$  be  sequences of real numbers such that for every $k\in {\mathbb{Z}}$, $ a_k<b_k$ and for every  $ j\in {\mathbb{Z}},$ $c_j<d_j.$  In this paper, we analyse the approximation properties of the following type of bivariate generalized Kantorovich sampling series: For $f\in C({\mathbb{R}}^2)$ (The class of all uniformly continuous and bounded functions on ${\mathbb{R}^2}$),
\begin{eqnarray*}
(K_w^{\varphi}f)(x,y)&=&\sum_{k=-\infty}^{\infty}\sum_{j=-\infty}^{\infty}\frac{w^2}{\Delta_{a_{k}}
\Delta_{c_{j}}}\varphi(wx-k,wy-j)\int_{\frac{k+a_k}{w}}^{\frac{k+b_k}{w}}
\int_{\frac{j+c_j}{w}}^{\frac{j+d_j}{w}}f(u,v)dvdu,
\end{eqnarray*}
where $\Delta_{a_{k}}=b_{k}-a_{k}$ and $\Delta_{c_{j}}=d_{j}-c_{j}.$

The purpose of this study is to obtain some approximation properties of
the Kantorovich type bivariate generalization of sampling operators $(K_w^{\varphi}f)_{w>0}.$ First, we give the pointwise convergence theorem and Voronovskaja type asymptotic theorem and then we obtain the degree of approximation and quantitative version of Voronovskaja type theorem of these operators. At the end we give some examples of kernels such as Box spline kernels and Bochner-Riesz kernel to which the theory can be applied.

\section{Preliminaries}

Let us denote by $C^{0}:=C^{0}({\mathbb{R}}^2)=C({\mathbb{R}}^2)$ the space of all uniformly continuous and bounded functions $f:{\mathbb{R}}^2\rightarrow {\mathbb{R}}$ with the usual supnorm $\|f\|_{\infty}$ and that by $C^{n}:=C^{n}({\mathbb{R}}^2)$ the space of all $n$-times continuously differentiable functions for which
\begin{eqnarray*}
\frac{\partial ^{|k|}f}{\partial x^{k_1}\partial y^{k_2}}\in C^{0},
\end{eqnarray*}
for $k=(k_1,k_2)$ with $|k|=k_1+k_2=n\geq 1$.

For $\delta>0$ and $(x,y)\in {\mathbb{R}}^2,$ $U_{\delta}(x,y)$ denotes the ball of radius $\delta$  centered at $(x,y).$ i.e., $U_{\delta}(x,y)=\{(u,v)\in   {\mathbb{R}}^2:(x-u)^2+ (y-v)^2\le \delta^{2}\}.$

Let $\varphi\in C^{0}$ be fixed. For any $\nu \in \mathbb{N}_{0}={\mathbb{N}}\cup\{0\},$ $h=(h_1,h_2)\in \mathbb{N}_{0}^2$ with $|h|=h_1+h_2=\nu,$ we define the algebraic moments as
\begin{eqnarray*}
m_{(h_1,h_2)}^{\nu}(\varphi):=\sum_{k=-\infty}^{\infty}\sum_{j=-\infty}^{\infty}
\varphi(x-k,y-j)(k-x)^{h_1}(j-y)^{h_2}
\end{eqnarray*}
and the absolute moments by
\begin{eqnarray*}
M_{(h_1,h_2)}^{\nu}(\varphi):=\sup_{(x,y)\in\mathbb{R}^2}\sum_{k=-\infty}^{\infty}\sum_{j=-\infty}^{\infty}|\varphi(x-k,y-j)|\, |k-x|^{h_1}\, |j-y|^{h_2}
\end{eqnarray*}
and
\begin{eqnarray*}
M_\nu(\varphi):=\max_{|h|=v}M_{(h_1,h_2)}^{\nu}(\varphi).
\end{eqnarray*}

Note that for $\mu, \nu\in\mathbb{N}{_0}$ with $\mu<\nu,$ $M_{\nu}(\varphi)<+\infty$ implies $M_{\mu}(\varphi)<+\infty.$ Indeed for $\mu<\nu,$ we have

\noindent
$\displaystyle\sum_{k=-\infty}^{\infty}\sum_{j=-\infty}^{\infty}|\varphi(x-k,y-j)||k-x|^{h_1}|j-y|^{h_2}$
\begin{eqnarray*}
&=&\sum\sum_{(k,j)\in U_1(x,y)}|\varphi(x-k,y-j)||k-x|^{h_1}|j-y|^{h_2}
\\&&+\sum\sum_{(k,j)\notin U_1(x,y)}|\varphi(x-k,y-j)||k-x|^{h_1}|j-y|^{h_2}=I_1+I_2.
\end{eqnarray*}
It is easy to see that $I_1\leq 4\|\varphi\|_{\infty}.$ Now, we obtain $I_2.$ Since $(k,j)\notin U_1(x,y),$ we have $|k-x|>1/2$ or $|j-y|>1/2.$ Now, we get
\begin{eqnarray*}
J_2&\leq& \sum\sum_{|k-x|>1/2}|\varphi(x-k,y-j)|\frac{|k-x|^{\nu-h_2}|j-y|^{h_2}}{|k-x|^{\nu-\mu}}\\&&+\sum\sum_{|j-y|>1/2}|\varphi(x-k,y-j)|\frac{|k-x|^{h_1}|j-y|^{\nu-h_1}}{|j-y|^{\nu-\mu}}\leq 2^{\nu-\mu+1}M_{\nu}(\varphi).
\end{eqnarray*}

When $\varphi$ has compact support, we immediately have that
$M_{\nu}(\varphi)<+\infty$ for every ${\nu}\in\mathbb{N}{_0}.$\\

We suppose that the following assumptions hold:
\begin{itemize}
\item[(i)] for every $(x,y)\in {\mathbb{R}}^2,$ we have
\begin{eqnarray*}
\sum_{k=-\infty}^{\infty}\sum_{j=-\infty}^{\infty}\varphi(x-k,y-j)=1,
\end{eqnarray*}
\item[(ii)] $M_2(\varphi)<+\infty$ and there holds
\begin{eqnarray*}
\lim_{R\rightarrow\infty}\sum\sum_{(k,j)\notin U_{R}(x,y)}|\varphi(x-k,y-j)|((k-x)^2+(j-y)^2))=0
\end{eqnarray*}
uniformly with respect to  $(x,y)\in{\mathbb{R}}^2,$
\item[(iii)]
for every $(x,y)\in{\mathbb{R}}^2,$
\begin{eqnarray*}
m_{(1,0)}^{1}(\varphi,x,y)&:=&\sum_{k=-\infty}^{\infty}\sum_{j=-\infty}^{\infty}
\varphi(x-k,y-j)(k-x)=0\\
m_{(0,1)}^{1}(\varphi,x,y)&:=&\sum_{k=-\infty}^{\infty}\sum_{j=-\infty}^{\infty}
\varphi(x-k,y-j)(j-y)=0\\
m_{(2,0)}^{2}(\varphi,x,y)&:=&\sum_{k=-\infty}^{\infty}\sum_{j=-\infty}^{\infty}
\varphi(x-k,y-j)(k-x)^2\\
m_{(0,2)}^{2}(\varphi,x,y)&:=&\sum_{k=-\infty}^{\infty}\sum_{j=-\infty}^{\infty}
\varphi(x-k,y-j)(j-y)^2\\
m_{(1,1)}^{2}(\varphi,x,y)&:=&\sum_{k=-\infty}^{\infty}\sum_{j=-\infty}^{\infty}
\varphi(x-k,y-j)(k-x)(j-y)
\end{eqnarray*}
\end{itemize}

\section{Main Results}
First, we give a pointwise convergence theorem at continuity points of the function $f$ for the bivariate operators.

\begin{thm}
For $f\in L^{\infty}({\mathbb{R}}^2),$
$\displaystyle{
\lim_{w\rightarrow\infty}(K_w^{\varphi}f)(x,y)=f(x,y)}$
at every point  $(x,y)$ of continuity of $f.$ Moreover, if the function is uniformly continuous and bounded on $\mathbb{R}{^2},$ then
$\displaystyle{
\lim_{w\rightarrow\infty}\|K_w^{\varphi}f-f\|=0.}$
\end{thm}
\begin{pf}
Let $\epsilon>0$ be fixed. By the continuity of $f$ at the point $(x,y),$ there exists $\delta>0$ such that $|f(u,v)-f(x,y)|\leq \epsilon,$ whenever $\sqrt{(u-x)^2+(v-y)^2}\leq \delta.$ \\

\noindent Now, we have

\noindent
$|(K_w^{\varphi}f)(x,y)-f(x,y)|$
\begin{eqnarray*}
&\leq&\sum\sum_{(\frac{k}{w},\frac{j}{w})\in U_{\frac{\delta}{2}}(x,y)}\frac{w^2}{\Delta_{a_{k}}
\Delta_{c_{j}}}|\varphi(wx-k,wy-j)|\int_{\frac{k+a_k}{w}}^{\frac{k+b_k}{w}}\int_{\frac{j+c_j}{w}}^{\frac{j+d_j}{w}}|f(u,v)-f(x,y)|dvdu\\&&+
\sum\sum_{(\frac{k}{w},\frac{j}{w})\notin U_{\frac{\delta}{2}}(x,y)}\frac{w^2}{\Delta_{a_{k}}
\Delta_{c_{j}}}|\varphi(wx-k,wy-j)|\int_{\frac{k+a_k}{w}}^{\frac{k+b_k}{w}}
\int_{\frac{j+c_j}{w}}^{\frac{j+d_j}{w}}|f(u,v)-f(x,y)|dvdu\\
&=&I_1+I_2.
\end{eqnarray*}
By the continuity of $f$ at the point $(x,y),$ we get
\begin{eqnarray*}
I_l\leq \epsilon\sum\sum_{(\frac{k}{w},\frac{j}{w})\in U_{\frac{\delta}{2}}(x,y)}|\varphi(wx-k,wy-j)|\leq \epsilon.  M_{0}(\varphi),
\end{eqnarray*}
for sufficiently large $w>0.$ Let $\delta>0$ be fixed. If $(\frac{k}{w},\frac{j}{w})\notin U_{\frac{\delta}{2}}(x,y)$ then
$\sqrt{(wx-k)^2+(wy-j)^2}>\frac{w\delta}{2}.$ Hence, we get
\begin{eqnarray*}
I_2\leq\frac{8\|f\|_{\infty}}{w^2\delta^2}\sum\sum_{(\frac{k}{w},\frac{j}{w})\notin U_{\delta}(x,y)}|\varphi(wx-k,wy-j)|(wx-k)^2+(wy-j)^2\rightarrow 0 \mbox{ as } w\rightarrow \infty
\end{eqnarray*}
and hence the first part of the theorem follows. The second part can be  proved similarly.
\end{pf}
\begin{rmk}Since functions in  $L^{\infty}({\mathbb{R}}^2)$ are locally Lebesgue integrable,  the generalized Kantorovich sampling series used in the above theorem are well defined for $f\in L^{\infty}({\mathbb{R}}^2).$
\end{rmk}

Next, we obtain a Voronovskaja type theorem for the bivariate Kantorovich sampling operators.
\begin{thm}\label{t1}
Let $f\in L^{\infty}({\mathbb{R}}^2)$ be such that $\frac{\partial f}{\partial x}, \frac{\partial f}{\partial y} $ exist at $(x,y)$ and let $\{a_k\}_{k\in {\mathbb{Z}}}, \{b_k\}_{k\in {\mathbb{Z}}},\{c_j\}_{j\in {\mathbb{Z}}}$ and $\{d_j\}_{j\in {\mathbb{Z}}}$ be bounded sequences such that $a_{k}+b_{k}=\alpha,$ $c_{j}+d_{j}=\beta,$   $b_k-a_k\geq l_{0}>0,$ $c_j-d_j\geq s_{0}>0$ and $\sup_{(k,j)}\{|a_{k}|,|b_{k}|,|c_{k}|,|d_{k}|\}\le K$ for some $K>0.$  Then, we have
\begin{eqnarray*}
\lim_{w\rightarrow\infty}w[(K_w^{\varphi}f)(x,y)-f(x,y)]=\frac{1}{2}\bigg(\alpha\frac{\partial f}{\partial x}(x,y)+\beta\frac{\partial f}{\partial y}(x,y)\bigg).
\end{eqnarray*}
\end{thm}
\begin{pf}
Let $\displaystyle K=\sup_{(k,j)}\{|a_k|,|b_k|,|c_j|,|d_j|\}.$ Let $(x,y)\in {\mathbb{R}}^2.$ By Taylor's theorem, we have
\begin{eqnarray*}
f(u,v)=f(x,y)+\frac{\partial f}{\partial x}(x,y)(u-x)+\frac{\partial f}{\partial y}(x,y)
(v-y)+h(u-x,v-y))((u-x)+(v-y)),
\end{eqnarray*}
for some bounded function $h$ such that $h(u,v)\rightarrow 0$ as $(u,v)\rightarrow (0,0).$\\

Thus, we have\\

\noindent
$(K_w^{\varphi}f)(x,y)-f(x,y)$
\begin{eqnarray*}
&=&\frac{\partial f}{\partial y}(x,y)\sum_{k=-\infty}^{\infty}\sum_{j=-\infty}^{\infty}\frac{w^2}{\Delta_{a_{k}}
\Delta_{c_{j}}}
\varphi(wx-k,wy-j)\int_{\frac{k+a_k}{w}}^{\frac{k+b_k}{w}}
\int_{\frac{j+c_j}{w}}^{\frac{j+d_j}{w}}(u-x)dvdu\\&&\nonumber+
\frac{\partial f}{\partial y}(x,y)\sum_{k=-\infty}^{\infty}\sum_{j=-\infty}^{\infty} \frac{w^2}{\Delta_{a_{k}}
\Delta_{c_{j}}}
\varphi(wx-k,wy-j)\int_{\frac{k+a_k}{w}}^{\frac{k+b_k}{w}}
\int_{\frac{j+c_j}{w}}^{\frac{j+d_j}{w}}(v-y)dvdu\\&&\nonumber+
\sum_{k=-\infty}^{\infty}\sum_{j=-\infty}^{\infty} \frac{w^2}{\Delta_{a_{k}}
\Delta_{c_{j}}}
\varphi(wx-k,wy-j)\\&&\times\int_{\frac{k+a_k}{w}}^{\frac{k+b_k}{w}}
\int_{\frac{j+c_j}{w}}^{\frac{j+d_j}{w}}h(u-x,v-y)((u-x)+(v-y))dvdu\\
&=&I_1+I_2+I_3, (\mbox{say}).
\end{eqnarray*}

First, we obtain $I_1.$
\begin{eqnarray*}
I_1&=&\frac{\partial f}{\partial x}(x,y)\sum_{k=-\infty}^{\infty}\sum_{j=-\infty}^{\infty}\frac{w^2}{\Delta_{a_{k}}
\Delta_{c_{j}}}
\varphi(wx-k,wy-j)\int_{\frac{k+a_k}{w}}^{\frac{k+b_k}{w}}
\int_{\frac{j+c_j}{w}}^{\frac{j+d_j}{w}}(u-x)dvdu\\
&=&\frac{\partial f}{\partial x}(x,y)\sum_{k=-\infty}^{\infty}\sum_{j=-\infty}^{\infty}\frac{w}{2\Delta_{a_{k}}}\varphi(wx-k,wy-j)
\bigg[\bigg(\frac{k+b_k}{w}-x\bigg)^2-\bigg(\frac{k+a_k}{w}-x\bigg)^2\bigg]\\
&=&\frac{\partial f}{\partial x}(x,y)\frac{1}{2w}\sum_{k=-\infty}^{\infty}\sum_{j=-\infty}^{\infty}\varphi(wx-k,wy-j)[(b_k+a_k)+2(k-wx)]\\
&=&\frac{\partial f}{\partial x}(x,y)\frac{\alpha}{2w}\sum_{k=-\infty}^{\infty}\sum_{j=-\infty}^{\infty}\varphi(wx-k,wy-j)\\&&+\frac{\partial f}{\partial x}(x,y)\frac{1}{w}\sum_{k=-\infty}^{\infty}\sum_{j=-\infty}^{\infty}\varphi(wx-k,wy-j)(k-wx)\\
&=& \frac{\partial f}{\partial x}(x,y)\frac{\alpha}{2w}+\frac{1}{w}\frac{\partial f}{\partial x}(x,y)m_{(1,0)}^{1}(\varphi,wx,wy)=\frac{\alpha}{2w} \frac{\partial f}{\partial x}(x,y).
\end{eqnarray*}
Similarly, we obtain $\displaystyle I_2=\frac{\beta}{2w}\frac{\partial f}{\partial y}(x,y).$

Now, we estimate $I_3.$ Let $\epsilon>0$ be fixed. Then, there exists $\delta>0$ such that $|h(t,s)|\leq \epsilon$ for every $\sqrt{t^2+s^2}\leq\delta. $ Now, we have
\begin{eqnarray*}
|I_3|&\leq&\sum\sum_{(\frac{k}{w},\frac{j}{w})\in U_{\frac{\delta}{2}}(x,y)}
\frac{w^2}{\Delta_{a_{k}}
\Delta_{c_{j}}}|\varphi(wx-k,wy-j)|\int_{\frac{k+a_k}{w}}^{\frac{k+b_k}{w}}
\int_{\frac{j+c_j}{w}}^{\frac{j+d_j}{w}}|h(u-x,v-y)|(|u-x|+|v-y|)dvdu\\&&+
\sum\sum_{(\frac{k}{w},\frac{j}{w})\notin U_{\frac{\delta}{2}}(x,y)}
\frac{w^2}{\Delta_{a_{k}}
\Delta_{c_{j}}}|\varphi(wx-k,wy-j)|\int_{\frac{k+a_k}{w}}^{\frac{k+b_k}{w}}
\int_{\frac{j+c_j}{w}}^{\frac{j+d_j}{w}}|h(u-x,v-y)|(|u-x|+|v-y|)dvdu\\
&=&I_3'+I_3''.
\end{eqnarray*}
First, we estimate $I_3'.$ For suitable large $w,$
\begin{eqnarray*}
|I_3'|&\leq&\epsilon\sum\sum_{(\frac{k}{w},\frac{j}{w})\in U_{\frac{\delta}{2}}(x,y)}
\frac{w^2}{\Delta_{a_{k}}
\Delta_{c_{j}}}|\varphi(wx-k,wy-j)|\int_{\frac{k+a_k}{w}}^{\frac{k+b_k}{w}}
\int_{\frac{j+c_j}{w}}^{\frac{j+d_j}{w}}|u-x|dvdu\\&&+\epsilon
\sum\sum_{(\frac{k}{w},\frac{j}{w})\in U_{\frac{\delta}{2}}(x,y)}
\frac{w^2}{\Delta_{a_{k}}
\Delta_{c_{j}}}|\varphi(wx-k,wy-j)|\int_{\frac{k+a_k}{w}}^{\frac{k+b_k}{w}}
\int_{\frac{j+c_j}{w}}^{\frac{j+d_j}{w}}|v-y|dvdu\\
&:=&J_1+J_2.
\end{eqnarray*}
\begin{eqnarray*}
J_1&\leq&\frac{\epsilon }{2}\sum\sum_{(\frac{k}{w},\frac{j}{w})\in U_{\frac{\delta}{2}}(x,y)}\frac{w}{\Delta_{a_{k}}}
|\varphi(wx-k,wy-j)|\bigg(\bigg(\frac{k+b_k}{w}-x\bigg)^2+\bigg(\frac{k+a_k}{w}-x\bigg)^2\bigg)\\
&\leq& \frac{\epsilon }{2wl_0}\sum\sum_{(\frac{k}{w},\frac{j}{w})\in U_{\frac{\delta}{2}}(x,y)}
|\varphi(wx-k,wy-j)|\\&&\times\bigg((a_k^2+b_k^2)+2(a_k+b_k)|k-wx|+2(k-wx)^2\bigg)\\
&\leq& \frac{\epsilon }{2wl_0}\bigg(2K^2M_{0}(\varphi)+2|\alpha| M_{(1,0)}^{1}(\varphi)+2M_{(2,0)}^{2}(\varphi)\bigg).
\end{eqnarray*}
Similarly, we get
\begin{eqnarray*}
J_2&\leq& \frac{\epsilon }{2ws_0}\bigg(2K^2M_{0}(\varphi)+2|\beta|M_{(0,1)}^{1}(\varphi)+2M_{(0,2)}^{2}(\varphi)\bigg).
\end{eqnarray*}
Therefore, we have
\begin{eqnarray*}
|I_3'|&\leq& \frac{\epsilon }{w}\bigg(K^2M_{0}(\varphi)\bigg(\frac{1}{l_0}+\frac{1}{s_0}\bigg)+\frac{\alpha M_{(1,0)}^{1}(\varphi)}{l_0}+\frac{\beta M_{(0,1)}^{1}(\varphi)}{s_0}\\&&+\frac{M_{(2,0)}^{2}(\varphi)}{l_0}+\frac{M_{(0,2)}^{2}(\varphi)}{s_0}\bigg).
\end{eqnarray*}
Let $M>0$ be a constant such that $|h(t,s)|\leq M.$ Let $\delta>0$ be fixed. If $(\frac{k}{w},\frac{j}{w})\notin U_{\frac{\delta}{2}}(x,y),$ then $\displaystyle\sqrt{\bigg(\frac{k}{w}-x\bigg)^2+\bigg(\frac{j}{w}-y\bigg)^2}>\frac{\delta}{2}$ and by assumption (ii), we have
\begin{eqnarray*}
|I_3''|
&\leq& 4M \sum\sum_{(\frac{k}{w},\frac{j}{w})\notin U_{\frac{\delta}{2}}(x,y)}|\varphi(wx-k,wy-j)|\left(\frac{((k-wx)^2+(j-wy)^2)}{\delta w^2} \right)
\end{eqnarray*}
for sufficiently large $w.$
Hence, we get $\displaystyle\lim_{w\rightarrow \infty}wI_3=0.$ So the assertion follows.
\end{pf}

\begin{rmk}
The boundedness assumption on $f$ in the previous theorem can be relaxed by assuming that there are two positive constants $a,b$ such that $|f(x,y)|\leq a+b(x^2+y^2),$ for every $(x,y)\in\mathbb{R}{^2}.$
\end{rmk}
\begin{pf}
One can easily verify that
\begin{eqnarray*}
|K_w^{\varphi}f)(x,y)|&\leq&\sum_{k=-\infty}^{\infty}\sum_{j=-\infty}^{\infty}\frac{w^2}{\Delta_{a_{k}}
\Delta_{c_{j}}}|\varphi(wx-k,wy-j)|\int_{\frac{k+a_k}{w}}^{\frac{k+b_k}{w}}
\int_{\frac{j+c_j}{w}}^{\frac{j+d_j}{w}}|f(u,v)|dvdu\\
&\leq&\sum_{k=-\infty}^{\infty}\sum_{j=-\infty}^{\infty}\frac{w^2}{\Delta_{a_{k}}
\Delta_{c_{j}}}|\varphi(wx-k,wy-j)|\int_{\frac{k+a_k}{w}}^{\frac{k+b_k}{w}}
\int_{\frac{j+c_j}{w}}^{\frac{j+d_j}{w}} (a+b(u^2+v^2))dvdu\\
&\leq&M_0(\varphi)\bigg(a+b(x^2+y^2)+\frac{b}{w}(|\alpha||x|+|\beta||y|)+\frac{2bK^2}{w^2}\bigg)
\\&&+\frac{2b}{w}\left(|x|M_{(1,0)}^{1}(\varphi)+|y|M_{(0,1)}^{1}(\varphi)\right)
+\frac{b}{w^2}\left(|\alpha|M_{(1,0)}^{1}(\varphi)+|\beta|M_{(0,1)}^{1}(\varphi)\right)\\&&+
\frac{b}{w^2}\left(M_{(2,0)}^{2}(\varphi)+M_{(0,2)}^{2}(\varphi)\right)
\end{eqnarray*}
and hence the series $K_w^{\varphi}f$ is absolutely convergent for every $(x,y)\in\mathbb{R}{^2}.$ Moreover, for a fixed $(x,y)\in\mathbb{R}{^2},$
\begin{eqnarray*}
P_1(u,v)=f(x,y)+\frac{\partial f}{\partial x}(x,y)(u-x)+\frac{\partial f}{\partial y}(x,y)(v-y),
\end{eqnarray*}
the Taylor's polynomial of first order centered at the point $(x,y)$ by the Taylor's formula, we can write
\begin{eqnarray*}
\frac{f(u,v)-P_1(u,v)}{|u-x|+|v-y|}=h(u-x,v-y),
\end{eqnarray*}
where $h$ is a function such that $\displaystyle\lim_{(u,v)\rightarrow 0}h(u,v)=0.$ Then $h$ is bounded on a neighbourhood of $(x,y),$ say $U_\delta(x,y)>0.$ For $(u,v)\notin U_\delta(x,y),$ we have
\begin{eqnarray*}
|h(u-x,v-y)|\leq \frac{|f(u,v)|}{|u-x|+|v-y|}+ \frac{|f(x,y)|}{|u-x|+|v-y|}+ \left|\frac{\partial f}{\partial x}\right| \frac{|u-x|}{|u-x|+|v-y|}  + \left|\frac{\partial f}{\partial y}\right| \frac{|v-y|}{|u-x|+|v-y|},
\end{eqnarray*}
and all the  terms on the right-hand side of the above inequality are bounded for $(u,v)\notin U_\delta.$ Hence, $h(.-x, .-y)$ is bounded on $\mathbb{R}{^2}.$ The same Voronovskaya formula can be obtained by writing the proof along the lines of the proof of Theorem \ref{t1}.
\end{pf}

\begin{rmk}
One can easily verify that the following set of sequences satisfy the conditions of the above theorem: $a_k=\frac{1}{|k|+1}, b_k=2-\frac{1}{|k|+1}, c_j=-\frac{1}{|j|+1}$ and $d_j=1+\frac{1}{|j|+1}.$ Moreover, several sequences can be constructed in a similar way.
\end{rmk}

For $f\in C(\mathbb{R}{^2}),$ the modulus of continuity for the bivariate case is defined
as follows:
\begin{eqnarray*}
\omega(f,\delta_1,\delta_2):=\sup\bigg\{|f(t,s)-f(x,y)|: (t,s),(x,y)\in {\mathbb{R}{^2}}\,\ \mbox{and}\,\ |t-x|\leq \delta_1,\,\ |s-y|\leq \delta_2\bigg\},
\end{eqnarray*}
where $\delta_1>0$ and $\delta_2>0.$ Further, $\omega(f,\delta_1, \delta_2)$ satisfies the following properties:
\begin{itemize}
\item[(a)] $\displaystyle\omega(f,\delta_1, \delta_2)\rightarrow 0 \mbox{ if}\,\ \delta_1\rightarrow 0 \,\ \mbox{and}\,\ \delta_2\rightarrow 0;$\\
\item[(b)] $\displaystyle |f(t,s)-f(x,y)|\leq \omega(f,\delta_1, \delta_2)\bigg(1+\frac{|t-x|}{\delta_1}\bigg)\bigg(1+\frac{|s-y|}{\delta_2}\bigg);$
\end{itemize}
The details of the modulus of continuity for the bivariate case can be found
in \cite{ana}.

In what follows, we shall define the following $K$-functional:
\begin{eqnarray*}
K(f,\delta,C^{0},C^{(1)}):=\displaystyle\inf\bigg\{\|f-g\|_{\infty}+\delta\max\bigg\{\bigg\|\frac{\partial g}{\partial x}\bigg\|_{\infty}, \bigg\|\frac{\partial g}{\partial y}\bigg\|_{\infty}\bigg\}, g\in C^{(1)}\bigg\},
\end{eqnarray*}
where $f\in C^{0}$ and $\delta>0.$ For every $f\in C^{0}$ there holds
\begin{eqnarray}\label{e3}
K(f,\delta/2,C^{0},C^{1})=\frac{1}{2}\overline{\omega}(f,\delta),
\end{eqnarray}
where $\overline{\omega}(f,.)$ denotes the least concave majorant of ${\omega}(f,.)$ (see e.g. \cite{ana}).\\

Now, we give the estimate of the rate of convergence of the bivariate operators $(K_w^{\varphi}f).$
\begin{thm}
Let $f\in C(\mathbb{R}{^2})$ and $\displaystyle{\sup_{k,j}\{a_k, b_k,c_{j},d_{j}\}\le K}$ for some $K>0$. Then, for all $(x,y)\in\mathbb{R}{^2},$ we have
\begin{eqnarray*}
|(K_w^{\varphi}f)(x,y)-f(x,y)|&\leq &\omega(f,\delta_1,\delta_2)
\bigg(M_{0}(\varphi)+\frac{K M_{0}(\varphi)}{w}\bigg(\frac{1}{\delta_{1}}+ \frac{1}{\delta_{2}}+ \frac{K}{w\delta_{1}\delta_{2}}\bigg)  \\ &&+\frac{K}{w}\left(M^{1}_{(1,0)}+M^{1}_{(0,1)}\right)\bigg(\frac{1}{\delta_{1}}+\frac{1}
{\delta_{2}}+\frac{1}{\delta_{1}\delta_{2} w}\bigg)\bigg)
\end{eqnarray*}
for any $w>0,$ and $\delta_1>0,$ $\delta_2>0.$
\end{thm}
\begin{pf}
Using the linearity and positivity of the operators and taking into account the property (b) above, we have

\noindent
$|(K_w^{\varphi}f)(x,y)-f(x,y)|$
\begin{eqnarray*}
&\leq&\sum_{k=-\infty}^{\infty}\sum_{j=-\infty}^{\infty}\frac{w^2}{\Delta_{a_{k}}
\Delta_{c_{j}}}|\varphi(wx-k,wy-j)|\int_{\frac{k+a_k}{w}}^{\frac{k+b_k}{w}}
\int_{\frac{j+c_j}{w}}^{\frac{j+d_j}{w}}|f(u,v)-f(x,y)|dvdu\\
&\leq& \omega(f,\delta_1, \delta_2)\sum_{k=-\infty}^{\infty}\sum_{j=-\infty}^{\infty}\frac{w^2}{\Delta_{a_{k}}
\Delta_{c_{j}}}|\varphi(wx-k,wy-j)|\\&&\times\int_{\frac{k+a_k}{w}}^{\frac{k+b_k}{w}}
\int_{\frac{j+c_j}{w}}^{\frac{j+d_j}{w}}
\bigg(1+\frac{|u-x|}{\delta_1}\bigg)\bigg(1+\frac{|v-y|}{\delta_2}\bigg)dvdu:=I_1+I_2+I_3+I_4.
\end{eqnarray*}
It is easy to see that $I_1=\omega(f,\delta_1,\delta_2)M_{0}(\varphi).$

Now, we estimate $I_2.$ For every $u\in\displaystyle\bigg[\frac{k+a_k}{w},\frac{k+b_k}{w}\bigg],$ we have
\begin{eqnarray*}
|u-x|&\leq&\bigg|u-\bigg(\frac{k+a_k}{w}\bigg)+\bigg(\frac{k+a_k}{w}\bigg)-x\bigg|
\leq \bigg|x-\frac{k}{w}\bigg|+ \frac{\max\{|a_{k}|,|b_{k}|\}}{w}\leq \bigg|x-\frac{k}{w}\bigg|+ \frac{K}{w} .
\end{eqnarray*}
Then, we obtain
\begin{eqnarray*}
I_2&\leq& \frac{\omega(f,\delta_1, \delta_2)}{\delta_1}\sum_{k=-\infty}^{\infty}\sum_{j=-\infty}^{\infty}\frac{w^2}{\Delta_{a_{k}}\Delta_{c_{j}}}
|\varphi(wx-k,wy-j)|\int_{\frac{k+a_k}{w}}^{\frac{k+b_k}{w}} \int_{\frac{j+c_j}{w}}^{\frac{j+d_j}{w}}|u-x|dvdu\\
&\leq& \frac{\omega(f,\delta_1, \delta_2)}{\delta_1 w}\sum_{k=-\infty}^{\infty}\sum_{j=-\infty}^{\infty}|\varphi(wx-k,wy-j)|\bigg(|k-wx|+K\bigg)\\
&\leq& \frac{\omega(f,\delta_1, \delta_2)}{\delta_1 w}\bigg(K M_{0}(\varphi)+M_{(0,1)}^{1}(\varphi)\bigg).
\end{eqnarray*}
Similarly, we obtain $I_3\le\displaystyle\frac{\omega(f,\delta_1, \delta_2)}{\delta_1 w}\bigg(K M_{(0,0)}^{0}(\varphi)+M_{(1,0)}^{1}(\varphi)\bigg).$\\
Now,
\begin{eqnarray*}
I_4&\leq&\frac{\omega(f,\delta_1, \delta_2)}{\delta_1 \delta_2}\frac{1}{w^2}
\sum_{k=-\infty}^{\infty}\sum_{j=-\infty}^{\infty}|\varphi(wx-k,wy-j)|\bigg(|k-wx|+K \bigg)\bigg(|j-wy|+K \bigg)\\
&\leq&\frac{\omega(f,\delta_1, \delta_2)}{\delta_1 \delta_2}\frac{1}{w^2}\bigg(K^2 M^{0}(\varphi) +K\left( M_{(1,0)}^{1}(\varphi) + M_{(0,1)}^{1}(\varphi)\right)+ M_{(1,1)}^{2}(\varphi)\bigg).
\end{eqnarray*}
Combining the estimates $I_1-I_4,$ we get the desired result.
\end{pf}

Using an estimate of the remainder in two dimensional Taylor formula and a technique developed in \cite{gonska} by Gonska et al., we obtain a quantitative version of the Voronovskaja formula for bivariate sampling series. We remark that quantitative
Voronovskaja formulae have important links with the theory of semi-groups
of operators (see \cite{FA1,Bern}). The strict connections between the two theories
were described in \cite{sem1} and recently developed in \cite{sem2,sem3,sem4}.

Now, we have the following quantitative version of the Voronovskaja formula for bivariate sampling series.

\begin{lem}\label{l1}
Let $f\in C^{1}$ and $(x,y)\in\mathbb{R}{^2},$ $(x_0,y_0)\in\mathbb{R}{^2}.$ Then for the remainder $R_1(f):=R_1(f,(x_0,y_0),(x,y))$ in the Taylor formula, we have
\begin{eqnarray*}
|R_1(f)|\leq |x-x_0|\overline{\omega}\bigg(\frac{\partial f}{\partial x}, |x-x_0|\bigg)+
|y-y_0|\overline{\omega}\bigg(\frac{\partial f}{\partial y}, |y-y_0|\bigg).
\end{eqnarray*}
\end{lem}

\begin{pf}
The proof follows from the Taylor formula of the first order using the same analysis as of Lemma 1 in \cite{bv1}.
\end{pf}

\begin{thm}\label{t2}
Let $f\in C^{1}(\mathbb{R}{^2})$ and let $a_{k},b_{k}, c_{j}, d_{j}$  satisfy the conditions of
 theorem \ref{t1}.  Then, for any $(x,y)\in\mathbb{R}{^2},$ we have

\noindent
$\displaystyle\left|w[(K_w^{\varphi}f)(x,y)-f(x,y)]-\frac{1}{2}\bigg(\alpha\frac{\partial f}{\partial x}+\beta\frac{\partial f}{\partial y}\bigg)\right|$
\begin{eqnarray*}
\leq A\overline{\omega}\bigg(\bigg\|\frac{\partial f}{\partial x}\bigg\|_{\infty}, \frac{B}{Aw}\bigg)+
C\overline{\omega}\bigg(\bigg\|\frac{\partial f}{\partial y}\bigg\|_{\infty}, \frac{D}{Cw}\bigg),
\end{eqnarray*}
where $A=\frac{1}{ l_0}
\bigg(2K^{2} M_{0}(\varphi)+2\alpha M_{(1,0)}^{1}(\varphi)+ 2 M_{(2,0)}^{2}(\varphi)\bigg),$
$B=
\bigg(K^{2}M_{0}(\varphi)+ \alpha M_{(1,0)}^{1}(\varphi)+M_{(2,0)}^{2}(\varphi)\bigg),$
$C=\frac{1}{ s_0}
\bigg(2K^{2} M_{0}(\varphi)+2\beta M_{(0,1)}^{1}(\varphi)+ 2 M_{(0,2)}^{2}(\varphi)\bigg)$ and
$D=
\bigg(K^{2}M_{0}(\varphi)+ \beta M_{(0,1)}^{1}(\varphi)+M_{(0,2)}^{2}(\varphi)\bigg).$
\end{thm}

\begin{pf}
Let $f\in C^{1}$ be fixed. Then, we can write

\noindent
$\left|w\left[(K_w^{\varphi}f)(x,y)-f(x,y)\right]-\frac{1}{2}\bigg(\alpha\frac{\partial f}{\partial x}+\beta\frac{\partial f}{\partial y}\bigg)\right|$
\begin{eqnarray*}
&=&\bigg|\frac{\partial f}{\partial x}(x,y)\sum_{k=-\infty}^{\infty}\sum_{j=-\infty}^{\infty}\frac{w^2}{\Delta_{a_{k}}\Delta_{c_{j}}}
\varphi(wx-k,wy-j)\int_{\frac{k+a_k}{w}}^{\frac{k+b_k}{w}} \int_{\frac{j+c_j}{w}}^{\frac{j+d_j}{w}}(u-x)dvdu\\&&
+\frac{\partial f}{\partial y}(x,y)\sum_{k=-\infty}^{\infty}\sum_{j=-\infty}^{\infty}\frac{w^2}{\Delta_{a_{k}}\Delta_{c_{j}}}\varphi(wx-k,wy-j)
\int_{\frac{k+a_k}{w}}^{\frac{k+b_k}{w}} \int_{\frac{j+c_j}{w}}^{\frac{j+d_j}{w}}(v-y)dvdu\\&&+
\sum_{k=-\infty}^{\infty}\sum_{j=-\infty}^{\infty}\frac{w^2}{\Delta_{a_{k}}\Delta_{c_{j}}}\varphi(wx-k,wy-j)\\&&\times
\int_{\frac{k+a_k}{w}}^{\frac{k+b_k}{w}} \int_{\frac{j+c_j}{w}}^{\frac{j+d_j}{w}}h(u-x,v-y)((u-x)+(v-y))dudv-\frac{1}{2}\bigg(\alpha\frac{\partial f}{\partial x}+\beta\frac{\partial f}{\partial y}\bigg)\bigg|\\
&\leq&\sum_{k=-\infty}^{\infty}\sum_{j=-\infty}^{\infty}\frac{w^2}{\Delta_{a_{k}}\Delta_{c_{j}}}|
\varphi(wx-k,wy-j)|\\&&\times\int_{\frac{k+a_k}{w}}^{\frac{k+b_k}{w}} \int_{\frac{j+c_j}{w}}^{\frac{j+d_j}{w}}|h(u-x,v-y)|(|u-x|+|v-y|)dvdu:=I_1.
\end{eqnarray*}
Using the relation (\ref{e3}) and Lemma \ref{l1}, we obtain
\begin{eqnarray*}
I_1&\le &\sum_{k=-\infty}^{\infty}\sum_{j=-\infty}^{\infty}\frac{w^2}{\Delta_{a_{k}}\Delta_{c_{j}}}|\varphi(wx-k,wy-j)|
\int_{\frac{k+a_k}{w}}^{\frac{k+b_k}{w}} \int_{\frac{j+c_j}{w}}^{\frac{j+d_j}{w}} |u-x|\overline{\omega}\bigg(\frac{\partial f}{\partial x},|u-x|\bigg)dvdu\\&&+
\sum_{k=-\infty}^{\infty}\sum_{j=-\infty}^{\infty}\frac{w^2}{\Delta_{a_{k}}\Delta_{c_{j}}}|\varphi(wx-k,wy-j)|
\int_{\frac{k+a_k}{w}}^{\frac{k+b_k}{w}} \int_{\frac{j+c_j}{w}}^{\frac{j+d_j}{w}}|v-y|\overline{\omega}\bigg(\frac{\partial f}{\partial y},|v-y|\bigg)dvdu\\
&:=&I_1'+I_1''.
\end{eqnarray*}
First, we estimate $I_1'.$ For $g\in \bf{C^2},$ we have
\begin{eqnarray*}
|I_1'|&\leq&2\sum_{k=-\infty}^{\infty}\sum_{j=-\infty}^{\infty}\frac{w^2}{\Delta_{a_{k}}\Delta_{c_{j}}}
|\varphi(wx-k,wy-j)|\int_{\frac{k+a_k}{w}}^{\frac{k+b_k}{w}} \int_{\frac{j+c_j}{w}}^{\frac{j+d_j}{w}}|u-x|K\bigg(\frac{\partial f}{\partial x},\frac{|u-x|}{2}\bigg)dvdu\\
&\leq&\sum_{k=-\infty}^{\infty}\sum_{j=-\infty}^{\infty}\frac{2w^2}{\Delta_{a_{k}}\Delta_{c_{j}}}
|\varphi(wx-k,wy-j)|\\&&\times \int_{\frac{k+a_k}{w}}^{\frac{k+b_k}{w}}\int_{\frac{j+c_j}{w}}^{\frac{j+d_j}{w}}|u-x|\bigg(\left\|\frac{\partial (f-g)}{\partial x}\right\|_{\infty}+\frac{|u-x|}{2}\max\bigg\{\bigg\|\frac{\partial g}{\partial x}\bigg\|_{\infty},\bigg\|\frac{\partial g}{\partial y}\bigg\|_{\infty}\bigg\}\bigg)dvdu\\
&\leq&\left\|\frac{\partial (f-g)}{\partial x}\right\|_{\infty}\sum_{k=-\infty}^{\infty}\sum_{j=-\infty}^{\infty}\frac{2w^2}{\Delta_{a_{k}}\Delta_{c_{j}}}
|\varphi(wx-k,wy-j)|\\&&\times \int_{\frac{k+a_k}{w}}^{\frac{k+b_k}{w}} \int_{\frac{j+c_j}{w}}^{\frac{j+d_j}{w}}|u-x|dvdu+
\sum_{k=-\infty}^{\infty}\sum_{j=-\infty}^{\infty}\frac{w^2}{\Delta_{a_{k}}\Delta_{c_{j}}}
|\varphi(wx-k,wy-j)|\\&&\times \int_{\frac{k+a_k}{w}}^{\frac{k+b_k}{w}} \int_{\frac{j+c_j}{w}}^{\frac{j+d_j}{w}} {(u-x)^2}\max\bigg\{\bigg\|\frac{\partial g}{\partial x}\bigg\|_{\infty},\bigg\|\frac{\partial g}{\partial y}\bigg\|_{\infty}\bigg\}dvdu:=J_1+J_2.
\end{eqnarray*}
It is easy to see that
\begin{eqnarray*}
|J_1|&\leq&\frac{w\left\|\frac{\partial (f-g)}{\partial x}\right\|_{\infty}}{\Delta_{a_{k}}}\sum_{k=-\infty}^{\infty}\sum_{j=-\infty}^{\infty}|\varphi(wx-k,wy-j)| \left|\bigg(\frac{k+b_k}{w}-x\bigg)^2+\bigg(\frac{k+a_k}{w}-x\bigg)^2\right|\\
&\leq&\left\|\frac{\partial (f-g)}{\partial x} \right\|_{\infty} \frac{1}{w l_0}
\bigg(2K^{2} M_{0}(\varphi)+2\alpha M_{(1,0)}^{1}(\varphi)+ 2 M_{(2,0)}^{2}(\varphi)\bigg)\\
\end{eqnarray*}
Now,
\begin{eqnarray*}
J_2&\leq&\max\left(\left\|\frac{\partial g}{\partial x}\right\|_{\infty}, \left\|\frac{\partial g}{\partial y}\right\|_{\infty}\right)\sum_{k=-\infty}^{\infty}\sum_{j=-\infty}^{\infty}\frac{w^2}{\Delta_{a_{k}}\Delta_{c_{j}}}|\varphi(wx-k,wy-j)| \int_{\frac{k+a_k}{w}}^{\frac{k+b_k}{w}} \int_{\frac{j+c_j}{w}}^{\frac{j+d_j}{w}}{(u-x)^2}dvdu\\
&\leq&\max\left(\left\|\frac{\partial g}{\partial x}\right\|_{\infty}, \left\|\frac{\partial g}{\partial y}\right\|_{\infty}\right)\frac{w}{3\Delta_{a_{k}}}\sum_{k=-\infty}^{\infty}\sum_{j=-\infty}^{\infty}|\varphi(wx-k,wy-j)| \left(\bigg(\frac{k+b_k}{w}-x\bigg)^3-\bigg(\frac{k+a_k}{w}-x\bigg)^3\right)\\
&\leq& \max\left(\left\|\frac{\partial g}{\partial x}\right\|_{\infty}, \left\|\frac{\partial g}{\partial y}\right\|_{\infty}\right)\frac{1}{w^2}
\bigg(K^{2}M_{0}(\varphi)+ \alpha M_{(1,0)}^{1}(\varphi)+M_{(2,0)}^{2}(\varphi)\bigg).
\end{eqnarray*}
Thus, we have
\begin{eqnarray*}
w I_1'\leq A\bigg(\left\|\frac{\partial (f-g)}{\partial x}\right\|_{\infty}+\frac{1}{w}\frac{B}{A}\max\left(\left\|\frac{\partial g}{\partial x}\right\|_{\infty}, \left\|\frac{\partial g}{\partial y}\right\|_{\infty}\right)\bigg).
\end{eqnarray*}
Taking the infimum over all $g\in C^2,$ we obtain
\begin{eqnarray*}
wI_1'\leq A\overline{\omega}\bigg(\left\|\frac{\partial (f-g)}{\partial x}\right\|_{\infty}, \frac{B}{A w}\bigg).
\end{eqnarray*}
Similarly, we obtain
\begin{eqnarray*}
wI_1''\leq C\overline{\omega}\bigg(\left\|\frac{\partial (f-g)}{\partial y}\right\|_{\infty}, \frac{D}{C w}\bigg).
\end{eqnarray*}
Combining the estimates $I_1'$ and $I_1'',$ we get the required result.
\end{pf}

\section{Applications to special kernels.}
In this section, we describe some particular examples of kernels $\varphi$ which illustrate the previous theory. In particular, we will examine the box splines kernel and Bochner-Riesz Kernel.
\subsection{Separable Box Splines}
First, we consider the sampling Kantorovich operators based on box spline functions which are separable. Let $\beta_{d_{1}}$ and $\beta_{d_{2}}$  be two Cardinal central B-splines of  degrees $d_{1}$ and $d_{2}$  respectively. i.e.,
\begin{eqnarray*}
\beta_{d_{1}}(x)&:=&\chi_{[-\frac{1}{2}, \frac{1}{2}]}\star \chi_{[-\frac{1}{2}, \frac{1}{2}]}\star \chi_{[-\frac{1}{2}, \frac{1}{2}]}\star\cdots \star\chi_{[-\frac{1}{2}, \frac{1}{2}]}\;\;\; (d_{1}+1 \,\ \mbox{times})\\
\beta_{d_{2}}(y)&:=&\chi_{[-\frac{1}{2}, \frac{1}{2}]}\star \chi_{[-\frac{1}{2}, \frac{1}{2}]}\star \chi_{[-\frac{1}{2}, \frac{1}{2}]}\star\cdots\star\chi_{[-\frac{1}{2}, \frac{1}{2}]}\;\;\;( d_{2}+1 \,\ \mbox{times}),
\end{eqnarray*}
where
\begin{eqnarray*}\chi_{[-\frac{1}{2}, \frac{1}{2}]}(x)=
\left\{
\begin{array}{ll}
1 , &\mbox{if} -\frac{1}{2}\leq x\leq\frac{1}{2}\\
0, & \mbox{otherwise}
\end{array}
\right.
\end{eqnarray*}
and $*$ denotes the convolution.\\
The Fourier transform of the functions $\beta_{d_{1}}(x)$  and $\beta_{d_{2}}(y)$ are given by
\begin{eqnarray*}
\widehat{\beta}_{d_{1}}(u)
=\bigg(\frac{\sin u/2}{u/2}\bigg)^{d_{1}+1}, \,\ u\in {\mathbb{R}},\\
\widehat{\beta}_{d_{2}}(v)
=\bigg(\frac{\sin v/2}{v/2}\bigg)^{d_{2}+1}, \,\ v\in {\mathbb{R}}.
\end{eqnarray*}
(see \cite{Stens2} and \cite{scho}). Given real numbers $\epsilon_0, \epsilon_1, \epsilon_{0}', \epsilon_1'$ with
$\epsilon_0<\epsilon_1,$  $\epsilon_0'<\epsilon_1'$ and  define $\varphi(x,y)= \varphi_{1}(x)\varphi_2(y),$ where
\begin{eqnarray*}
 \varphi_{1}(x)&=&a_0\beta_{d_{1}}(x-\epsilon_0)+a_1\beta_{d_{1}}(x-\epsilon_1) \\
 \varphi_2(y)&=&
 b_0\beta_{d_{2}}(y-\epsilon_0')+b_1\beta_{d_{2}}(y-\epsilon_1')
\end{eqnarray*}
The Fourier transform of $\varphi$ is   $\widehat{\varphi}(u,v)=\widehat{\varphi_{1}}(u) \widehat{\varphi_{2}}(v),$ where
\begin{eqnarray*}
\widehat{\varphi_{1}}(u)=\bigg(a_0e^{-i\epsilon_0 u}+a_1e^{-i\epsilon_1 u}\bigg)\widehat{\beta}_{d_{1}}(u)\mbox{ and }
\widehat{\varphi_{2}}(v)=\bigg(b_0e^{-i\epsilon_0' v}+b_1e^{-i\epsilon_1' v}\bigg)\widehat{\beta}_{d_{2}}(v).
\end{eqnarray*}
Using the Poisson summation formula
\begin{eqnarray*}
(-i)^{j}\sum_{k=-\infty}^{\infty}\varphi_{1}(u-k)(u-k)^{j}\sim \sum_{k=-\infty}^{\infty}\widehat{\varphi_{1}}^{(j)}(2\pi k)e^{i2\pi ku},
\end{eqnarray*}
we obtain
\begin{eqnarray*}
\sum_{k=-\infty}^{\infty}\varphi_{1}(u-k)=\sum_{k=-\infty}^{\infty}\widehat{\varphi_{1}}(2\pi k)e^{i2\pi ku}.
\end{eqnarray*}
We have
\begin{eqnarray*}
\widehat{\beta_{d_{1}}}(2\pi k)=\bigg(\frac{\sin(\pi k)}{\pi k}\bigg)^{d_{1}+1}
=
\left\{
\begin{array}{ll}
1 , &\mbox{if}\,\ k=0\\
0, & \mbox{if}\,\ k\neq 0
\end{array}
\right.
\end{eqnarray*}
and hence
\begin{eqnarray*}
\widehat{\varphi_{1}}(2\pi k)
=
\left\{
\begin{array}{ll}
a_0+a_1 , &\mbox{if}\,\ k=0\\
0, & \mbox{if}\,\ k\neq 0
\end{array}
\right.
\end{eqnarray*}
Thus
\begin{eqnarray*}
\sum_{k=-\infty}^{\infty}\varphi_{1}(u-k)=a_0+a_1.
\end{eqnarray*}
Therefore in order to satisfy condition (i), we assume  $a_0+a_1=1$ and $b_0+b_1=1.$ Now, we show that condition (iii) is also satisfied.

Again from the Poisson summation formula, we obtain
\begin{eqnarray*}
(-i)\sum_{k=-\infty}^{\infty}\varphi_{1}(u-k)(u-k)=
\sum_{k=-\infty}^{\infty}\widehat{\varphi_{1}}^{\prime}(2\pi k)e^{i2\pi ku}.
\end{eqnarray*}
Also, we have
\begin{eqnarray*}
\widehat{\varphi_{1}}^{\prime}(u)=(-i\epsilon_0a_0e^{-i\epsilon_0 u}-i\epsilon_1a_1e^{-i\epsilon_1 u})\widehat{\beta}_{d_{1}}(u)+(a_0e^{-i\epsilon_0 u}+a_1e^{-i\epsilon_1 u})\widehat{\beta}_{d_{1}}^{\prime}(u).
\end{eqnarray*}
Since $\widehat{\beta}_{d_{1}}^{\prime}(2\pi k)=0,$ $\forall k$ which implies that $\widehat{\varphi_{1}}^{\prime}(2\pi k)=0.$ Thus, we have
\begin{eqnarray*}
\widehat{\varphi_{1}}(0)=a_0+a_1=1,\,\,\ \widehat{\varphi_{1}}^{\prime}(0)=\epsilon_0a_0+\epsilon_1a_1=0.
\end{eqnarray*}
Solving the above linear system we get the unique solution
\begin{eqnarray*}
a_0=\frac{\epsilon_1}{\epsilon_1-\epsilon_0}, a_1=-\frac{\epsilon_0}{\epsilon_0-\epsilon_1}.
\end{eqnarray*}
Similarly we obtain
\begin{eqnarray*}
b_0=\frac{\epsilon_1'}{\epsilon_1'-\epsilon_0'}, b_1=-\frac{\epsilon_0'}{\epsilon_0'-\epsilon_1'}.
\end{eqnarray*}
It is easy to check that with these values of the parameters, the kernel function $\varphi$ satisfies all the conditions. Hence for these kernels, we can obtain the Voronovskaja formula as in Theorem 3. 1 and the corresponding quantitative estimate.

\subsection{Translates of Box Splines}
Boor and Devore \cite{Boor}, initiated the theory of bivariate splines and extensively studied by several researchers (see, e.g., \cite{Chui} and \cite{Hollig}). First, we give some basic definitions here.

Let $A=(a_{\nu,\mu})$ be a  $2\times m$ matrix with $m\geq 3$ and column vectors $A_{\mu}\in\mathbb{Z}{^2}\backslash\{(0,0)\},$ $\mu=1,2,\ldots ,m$ be such that rank(A)=2.
The box spline $B_{A}$ is defined  via the following:
\begin{eqnarray*}
\int_{\mathbb{R}{^2}}B_{A}(x,y)g(x,y)dxdy=\int_{Q^{m}}g(At)dt \;\;\;(g\in C({\mathbb{R}}^2)),
\end{eqnarray*}
where $Q^{m}=[-\frac{1}{2},\frac{1}{2}]^{m}$ and $t=(t_1,t_2...,t_m)\in Q^{m}.$ It is well known that $B_{A}(x,y)\geq 0,$ $(x,y)\in\mathbb{R}{^2}$ and $Supp(B_{A})=AQ^{m}$ (see \cite{Chui}). Moreover, if $\rho=\rho(A)$ is the greatest integer such that all submatrices generated from $A$ by deleting $\rho$ columns, have rank $2,$ then $B_{A}\in C^{(\rho-1)}(\mathbb{R}{^2})$ and the Fourier transform of $B_{A}$ is given by
\begin{eqnarray*}
\widehat{{B}}_{A}(u,v)=\frac{1}{2\pi}\prod_{\mu=1}^{m}sinc\bigg(\frac{1}{2\pi}(ua_{1\mu}+va_{2\mu})\bigg),\,\,\ (u,v)\in\mathbb{R}{^2}
\end{eqnarray*}
where $sinc(z)=\frac{\sin \pi z}{\pi z}$ for $z\in\mathbb{R}\backslash\{0\},$ $sinc(0)=1.$

The functions $B_{A}$ are piecewise polynomials of total degree $m-2.$ In order to construct kernels satisfying the conditions (i)-(iii) by using these box splines. we  have the following result (see \cite{PLB1} and \cite{Fis}).

\begin{prop}
Let $r,m\in\mathbb{N},$ $r\geq 2,$ $m\geq r+1$ and let $A$ be a $2\times m$ matrix with column vectors $A_{\mu}\in \mathbb{Z}{^2}\backslash \{(0,0)\}$ and $\rho(A)\geq r-1.$ If $b_{hr,}h=(h_1,h_2)\in \mathbb{N}{_0}^2$ with $|h|=h_1+h_2\leq r_0:=2[\frac{r-1}{2}]$ are the unique solutions of the linear system
\begin{eqnarray*}
(-1)^{\lfloor|s|/2\rfloor}\sum_{|h|\leq r_0}b_{hr}h_1^{s_1}h_2^{s_2}=\frac{\partial ^{|s|}}{\partial u^{s_1}\partial v^{s_2}}\frac{1}{\widehat{{B}}_{A}}(0,0),
\end{eqnarray*}
where $s=(s_1,s_2)\in \mathbb{N}{_0}^2 ,$ $|s|=s_1+s_2\leq r_0$ and $[\alpha]$ denotes the largest integer but not bigger than $\alpha,$ then
\begin{eqnarray*}
\varphi_{A,r}(x,y)=\frac{1}{2\pi}\bigg(b_{0r}B_{A}(x,y)+\frac{1}{2}\sum_{0<|h|\leq r_0}b_{hr}(B_{A}(x+h_1,y+h_2)+B_{A}(x-h_1,y-h_2))\bigg)
\end{eqnarray*}
is a polynomial spline of degree $m-2$ with compact support, $(r-2)$ times continuously differentiable satisfying the following assumptions:
\begin{itemize}
\item[(a)] $\displaystyle\sum_{k=-\infty}^{\infty} \sum_{j=-\infty}^{\infty}\varphi_{A,r}(x,y)(x-k,y-j)=1,$\,\,\ $(x,y)\in \mathbb{R}{^2}$\\
\item[(b)] $\displaystyle B_{r}(\varphi_{A,r})=\max_{|s|=r}B_{(s_1,s_2)}^{r}(\varphi_{A,r})<+\infty$\\
\item[(c)] $\displaystyle\sum_{k=-\infty}^{\infty} \sum_{j=-\infty}^{\infty}\varphi_{A,r}(x-k,y-j)(x-k)^{s_1}(y-j)^{s_2}=0,$ \,\,\ for all $(x,y)\in \mathbb{R}{^2},$ $s=(s_1,s_2)\in \mathbb{N}{_0}^2$ with $0<|s|<r.$
\end{itemize}
\label{prop4.1}
\end{prop}
The condition (c) of Proposition \ref{prop4.1} is equivalent to the following one which uses the Fourier transform of $\varphi_{A,r}.$ For every $s=(s_1,s_2)$ with $0\leq |s|<r,$
\begin{eqnarray*}
\frac{\partial ^{|s|}\widehat{\varphi}_{A,r}}{\partial u^{s_1}\partial v^{s_2}}(2k\pi,2j\pi)
=
\left\{
\begin{array}{ll}
0 , &\mbox{if}\,\ (k,j)\neq (0,0)\\
\delta_{s0}, & \mbox{if}\,\ (k,j)= (0,0),
\end{array}
\right.
\end{eqnarray*}
where
\begin{eqnarray*}
\delta_{s0}
=
\left\{
\begin{array}{ll}
1 , &\mbox{if}\,\ (s_1,s_2)=(0,0)\\
0, & \,\ \mbox{otherwise} .
\end{array}
\right.
\end{eqnarray*}
Now, we construct the linear combinations of translates of $\varphi_{A,r}$ for $r\geq 3$ in such a way that all the assumptions are satisfied for given values of the constants $S_1,$ $S_2$ and $S_3.$
Let $S_1,S_2$ and $S_3$ be real numbers such that $(S_1,S_2,S_3)\neq (0,0,0).$ For given constants $a,b,c,$ we define the kernel function $\varphi$  by
\begin{eqnarray*}
\varphi(x,y)&=&C_1\varphi_{A,r}(x,y)+C_2(\varphi_{A,r}(x+a,y+a)+\varphi_{A,r}(x-a,y-a))
\\&&+C_3(\varphi_{A,r}(x+b,y)+\varphi_{A,r}(x-b,y))+C_4(\varphi_{A,r}(x,y+c)+\varphi_{A,r}(x,y-c)).
\end{eqnarray*}
The Fourier transform  $\widehat{\varphi}(u,v)$ of $\varphi$ can be written as
\begin{eqnarray*}
\hat{\varphi}(u,v)=\widehat{\varphi_{A,r}}(u,v)(C_1+2C_2 \cos(a(u+v))+2C_3\cos(bu)+2C_4\cos(cv)).
\end{eqnarray*}
Thus we have the following linear system of equations
\begin{eqnarray*}
\hat{\varphi}(0,0)&=&C_1+2C_2+2C_3+2C_4=1\\
\frac{\partial\hat{\varphi}}{\partial u}(0,0)&=&\frac{\partial\hat{\varphi}}{\partial v}(0,0)=0\\
\frac{\partial^2\hat{\varphi}}{\partial u^2}(0,0)&=&-2(a^2C_2+b^2C_3)=S_1\\
\frac{\partial^2\hat{\varphi}}{\partial v^2}(0,0)&=&-2(a^2C_2+c^2C_4)=S_2\\
\frac{\partial^2\hat{\varphi}}{\partial u\partial v}(0,0)&=&-2a^2C_2=S_3
\end{eqnarray*}
while all the above derivatives are zero at every point $(2k\pi,2j\pi)$ with $(k,j)\neq (0,0).$
When $\varphi$ has compact support, the limit relation in assumption
ii) is automatically satisfied. Hence, all the conditions are satisfied. Further, the bivariate generalized Kantorovich sampling series now takes the form
\begin{eqnarray}\label{ex1}
(K_w^{\varphi}f)(x,y)&=&\sum_{k=-\infty}^{\infty}\sum_{j=-\infty}^{\infty}\frac{w^2}{\Delta_{a_{k}}\Delta_{c_{j}}}
\varphi_{A,r}(wx-k,wy-j)\nonumber\\&&\times \int_{\frac{k+a_k}{w}}^{\frac{k+b_k}{w}} \int_{\frac{j+c_j}{w}}^{\frac{j+d_j}{w}}f(u,v)dvdu,
\end{eqnarray}
for every $w>0.$ For the generalized sampling series (\ref{ex1}), we can obtain the Voronovskaja formula of Theorem 3.1 and the corresponding quantitative estimate.

\subsection{Radial Kernels: The Bochner-Riesz Kernel}
Next, we consider the two dimensional Bochner-Riesz Kernel $b^{\gamma}(x,y)$(see, e.g., \cite{PLB3,PLB2,PLB1,Nesselr,EMS}). It is defined by
\begin{eqnarray*}
\varphi(x,y)\equiv b^{\gamma}(x,y)=\frac{2^{\gamma}\Gamma(\gamma+1)}{2\pi}\left(\sqrt{x^2+y^2}\right)^{-1-\gamma}J_{1+\gamma}\left(\sqrt{x^2+y^2}\right)
\end{eqnarray*}
for $\gamma>1/2,$ where $J_{\lambda}$ is the Bessel function of order $\lambda.$ It is well known that
\begin{eqnarray*}
\sum_{k=-\infty}^{\infty} \sum_{j=-\infty}^{\infty}b^{\gamma}(x-k,y-j)=1, (x,y)\in \mathbb{R}{^2}
\end{eqnarray*}
and further  if $\gamma>5/2,$ then
\begin{eqnarray*}
\sup_{(x,y)\in \mathbb{R}{^2}}\sum_{k=-\infty}^{\infty} \sum_{j=-\infty}^{\infty}|b^{\gamma}(x-k,y-j)|((x-k)^2+(y-j)^2)<+\infty
\end{eqnarray*}
and so $M_{2}(b^{\gamma})<+\infty.$
The following proposition shows that the kernel function $b^{\gamma}(x,y)$ satisfies condition (ii).
\begin{prop}
\cite{Panbar}. For the function $b^{\gamma},$ $\gamma>5/2$ there holds
\begin{eqnarray*}
\lim_{R\rightarrow +\infty}\sum\sum_{(k,j)\notin U_{R}(x,y)}|b^{\gamma}(x-k,y-j)|((x-k)^2+(y-j)^2)=0,
\end{eqnarray*}
uniformly with respect to $(x,y)\in \mathbb{R}{^2}.$
\end{prop}
As a consequences of the following results proved in \cite{PLB2}, Lemma 3.2 and Lemma 2 in \cite{bv1}, we obtain that $b^{\gamma}$ satisfies condition (iii).
\begin{prop}
Let $\varphi: \mathbb{R}{^2}\rightarrow \mathbb{R}$ be a continuous and bounded function such that $M_{r}(\varphi)<+\infty$ for some $r\in\mathbb{N}{_0},$ $h=(h_1,h_2)\in\mathbb{N}{_0}^2$ with $|h|=h_1+h_2\leq r,$ the double series in the definition of $M_r(\varphi)$ being uniformly convergent on each compact subset of $\mathbb{R}{^2}$ for each $h$ such that $|h|=r.$ The following two assertions are equivalent for $c\in\mathbb{R}:$
\begin{itemize}
\item[(a)] $\displaystyle\sum_{k=-\infty}^{\infty} \sum_{j=-\infty}^{\infty}\varphi(x-k,y-j)(k-x)^{h_1}(y-j)^{h_2}=c,$ $\forall (x,y)\in \mathbb{R}{^2}$\\
\item[(b)] $\displaystyle\frac{\partial ^{|h|}\hat{\varphi}}{\partial u^{h_1}\partial v^{h_2}}(2k\pi,2j\pi)
=
\left\{
\begin{array}{ll}
(-1)^{|h|}\frac{c}{2\pi} , &\mbox{if}\,\ (k,j)=(0,0)\\
0, & \mbox{if}\,\ (k,j)\neq (0,0).
\end{array}
\right.$
\end{itemize}
\end{prop}

\begin{prop}
\cite{bv1}. For the function $b^{\gamma},$ $\gamma>5/2$ we have
\begin{eqnarray*}
m_{(1,0)}^1(b^{\gamma},x,y)&=&m_{(0,1)}^1(b^{\gamma},x,y)=0\\
m_{(2,0)}^2(b^{\gamma},x,y)&=&m_{(0,2)}^2(b^{\gamma},x,y)=-\frac{\gamma}{\pi}\\
m_{(1,1)}^2(b^{\gamma},x,y)&=&0.
\end{eqnarray*}
\end{prop}
Thus, Bochner-Riesz Kernel  $b^{\gamma}(x,y)$  satisfies all required the conditions.  Moreover, the bivariate generalized Kantorovich sampling series now takes the form
\begin{eqnarray}\label{ex2}
(K_w^{\varphi}f)(x,y)&=&\frac{2^{\gamma}\Gamma(\gamma+1)}{2\pi}\sum_{k=-\infty}^{\infty}\sum_{j=-\infty}^{\infty} \frac{w^2}{\Delta_{a_{k}}\Delta_{c_{j}}}
\left(\sqrt{(wx-k)^2+(wy-j)^2}\right)^{-1-\gamma}\nonumber\\&&\times J_{1+\gamma}\left(\sqrt{(wx-k)^2+(wy-j)^2}\right)
\int_{\frac{k+a_k}{w}}^{\frac{k+b_k}{w}} \int_{\frac{j+c_j}{w}}^{\frac{j+d_j}{w}}f(u,v)dvdu,
\end{eqnarray}
for every $w>0.$ For the generalized sampling series (\ref{ex2}), we can obtain the Voronovskaja formula as in Theorem 3. 1 and the corresponding quantitative estimate.

\end{document}